\documentclass[10pt]{article}

\usepackage[english]{babel}
\usepackage{amsmath}
\usepackage{amssymb}
\usepackage{bm}
\usepackage{bbm}
\usepackage{mathrsfs,color}
\usepackage{graphics,graphicx,theorem}
\usepackage{dsfont}
\usepackage{setspace}

\usepackage{mathtools}

\usepackage{natbib}
\usepackage{multirow}
\usepackage{hhline}
\usepackage{hyperref}

\hypersetup{
	colorlinks=true,
	urlcolor=blue,
	citecolor=blue,
	linktoc=all,
	linkcolor=red} 

\usepackage[scriptsize]{subfigure}
\allowdisplaybreaks

\makeatletter
\renewcommand\@biblabel[1]{}
\makeatother

\def\bSig\mathbf{\Sigma}

\newcommand{\qed}{$\square$}

\usepackage[top=1.2in,bottom=1.2in,left=1.2in,right=1.2in]{geometry}

\begin{document}

\title{\bf {\Large{A characterization of product-form exchangeable feature probability functions}}}
\author{
Marco Battiston\\
\normalsize{Bocconi University, Italy}\\
\normalsize{email: \texttt{marco.battiston@phd.unibocconi.it}}
\bigskip\\
Stefano Favaro\\
\normalsize{University of Torino and Collegio Carlo Alberto, Italy}\\
\normalsize{email: \texttt{stefano.favaro@unito.it}}
\bigskip\\
Daniel M. Roy\\
\normalsize{University of Toronto, Canada}\\
\normalsize{email: \texttt{droy@utstat.toronto.edu}}
\bigskip\\
Yee Whye Teh\\
\normalsize{University of Oxford, UK}\\
\normalsize{email: \texttt{y.w.teh@stats.ox.ac.uk}}
}
\date{}
\maketitle
\thispagestyle{empty}

\setcounter{page}{1}
\begin{center}
\textbf{Abstract} 
\end{center}
We characterize the class of exchangeable feature allocations assigning probability $V_{n,k}\prod_{l=1}^{k}W_{m_{l}}U_{n-m_{l}}$ to a feature allocation of $n$ individuals, displaying $k$ features with counts $(m_{1},\ldots,m_{k})$ for these features. Each element of this class is parametrized by a countable matrix $V$ and two sequences $U$ and $W$ of non-negative weights. Moreover, a consistency condition is imposed to guarantee that the distribution for feature allocations of $n-1$ individuals is recovered from that of $n$ individuals, when the last individual is integrated out. In Theorem 1.1, we prove that the only members of this class satisfying the consistency condition are mixtures of the Indian Buffet Process over its mass parameter $\gamma$ and mixtures of the Beta--Bernoulli model over its dimensionality parameter $N$. Hence, we provide a characterization of these two models as the only, up to randomization of the parameters, consistent exchangeable feature allocations having the required product form. \vspace*{.2in}

\noindent\textsc{Keywords}: {Exchangeable Feature Allocations, Gibbs-type partitions, Indian Buffet process.} 

\maketitle

\section{Introduction} \label{s:intro}

\emph{Feature allocations} are popular models in machine learning. In these models, we consider a set of $n$ individuals, each displaying a (possibly empty) set of features. Specifically, let $(\mathcal{X},\mathcal{B})$ be measurable space, representing the collection of all possible features. Each individual is described by a random finite subset $X_{i}$ of $\mathcal{X}$, collecting his features. Each feature $x \in \mathcal{X}$ can be shared by many individuals. Given a set of $n$ individuals, a \emph{feature allocation} describes the sharing of features among these individuals. A way of describing this sharing is to associate to each of the $k$ points in $\cup_{1 \leq i \leq n} X_{i}$ a subset of $[n]:=\{ 1, \ldots , n \}$, summarizing the individuals having that particular feature. We denote by $(f_{n,1},\ldots,f_{n,k})$ the subsets of $[n]$ representing each of the $k$ features and by $(m_{1},\ldots,m_{k})$ the cardinalities of these sets. \par 
A feature allocation is \emph{exchangeable} when its distribution is invariant under permutation of the indexes of the individual, i.e. the feature allocation induced by the random sets $(X_{1},\ldots,X_{n})$ is equal in distribution to that induced by $(X_{\sigma(1)},\ldots,X_{\sigma(n)})$, for all permutation $\sigma$ of $[n]$.  Moreover, as pointed out in \citet{Bro13}, it is usually convenient to assign an order to the $k$ features present in a feature allocation of $n$ individuals. A way of achieving this purpose is drawing $k$ values from a continuous distribution and ordering the $k$ features accordingly. The resulting feature allocation is said to be a \emph{randomly ordered feature allocation}. In \citet{Bro13}, the authors study the class of exchangeable randomly ordered feature allocations admitting as a sufficient statistics the vector $(m_{1},\ldots,m_{k})$, i.e., the class of randomly ordered exchangeable feature allocations of the form
\begin{equation}
P(f_{n,1},\ldots,f_{n,k})=\pi _{n}(m_{1},\ldots,m_{k}),
\end{equation}
for a symmetric function $\pi_{n}$, called an \emph{exchangeable feature probability punction} (EFPF), \citet{Bro13}. \par

When dealing with random exchangeable feature allocations, we also require consistency conditions that guarantee that the distribution of a feature allocation of $n$ individuals coincides with that of $n-1$ individuals, when the last individual is integrated out. When considering randomly ordered exchangeable feature allocations with EFPF, this consistency notion specializes to the condition 
\begin{equation} \label{consistency.form}
\begin{split}
\pi_{n} \left(m_{1},\ldots,m_{k}\right)=
\sum_{j=0}^{\infty}\binom{k+j}{j}\sum_{\underline{z}\in\left\{ 0,1\right\} ^{k}}\pi_{n+1}(m_{1}+z_{1},\ldots,m_{k}+z_{k},\underbrace{1,\ldots,1}_{j}).
\end{split}
\end{equation} 
Feature allocations satisfying this condition are said to be \emph{consisent}. \par
The most remarkable example of exchangeable consistent feature allocation with EFPF is the \emph{Indian Buffet Process} (IBP), initially introduced in \citet{Gri06}, in its one parameter version, and then extended to its two, \citet{Gha07}, and three parameters versions, \citet{Teh10}.  The EFPF of a 3-parameter $(\gamma,\theta,\alpha)$ IBP has the following form 
\begin{equation*}
\frac{1}{k!}\Big(\frac{\gamma}{(\theta+1)_{n-1\uparrow}}\Big)^{k}\exp\Big(-\sum_{i=1}^{n}\gamma\frac{\left(\alpha+\theta\right)_{i-1\uparrow}}{\left(1+\theta\right)_{i-1\uparrow}}\Big)\prod_{l=1}^{k}\left(1-\alpha\right)_{m_{l}-1\uparrow}\left(\theta+\alpha\right)_{n-m_{l}\uparrow},
\end{equation*}
where $(x)_{m\uparrow}$ denotes the rising factorial and the parameters must satisfy the conditions $\gamma>0$, $0\leq\alpha<1$, and $\theta>-\alpha$. The 2-parameter IBP is recovered when $\alpha$ is set equal to zero, and the 1 parameter IBP when we also impose $\theta=1$. For a review of the IBP and its applications in machine learning, the reader is refered to \citet{Gri11}. \par
The IBP is derived as the limit of a Beta--Bernoulli model, in \citet{Gri06}. This latter model is the counterpart of the IBP when the set of all possible features $\mathcal{X}$ has finite cardinality, $N$. The EFPF of a Beta--Bernoulli model with parameters $(N,\alpha,\theta)$ is
\begin{equation} \label{Beta--Bernoulli}
 \binom{N}{k} \left(\frac{-\alpha}{(\theta+\alpha)_{n\uparrow}}\right) ^{k} \left( \frac{(\theta+\alpha)_{n\uparrow} }{(\theta)_{n\uparrow} } \right) ^{N}  \prod_{i=1}^{k} (1-\alpha)_{m_{i}-1\uparrow} (\theta+\alpha)_{n-m_{i}\uparrow},
\end{equation}  
where $\alpha<0$ and $\theta>-\alpha$. In Appendix \ref{App.BetaBernoulli}, we provide a brief description of the Beta--Bernoulli model and a derivation of its EFPF. \par
Feature allocations are generalizations of partitions. Indeed, a random partition is the particular case of a random feature allocation in which each random set $X_{i}$ is a singleton with probability one. All notions just introduced (consistent, exchangeable, ordered feature allocation and EFPF) were first introduced for partitions and only recently extended to the feature allocation case. The reader is referred to \citet{Pit06} for a complete review of exchangeable random partitions. One the most important distribution for random partitions is the Ewens--Pitman formula, which is generalisation of the famous Ewens formula. Starting from this distribution, \citet{Gne06} considers a larger class of random partitions, having an exchangeable partition probability function (see \citet{Pit95} for a definition) with the same product form as the Ewens--Pitman formula, but allowing a more general parametrization, depending on a triangular array and on a sequence of non-negative weights. Theorem 12 of \citet{Gne06} characterizes all elements of this class of distributions for random partitions satisfying a consistency condition similar to \eqref{consistency.form}. The resulting class of distributions is termed \emph{Gibbs-type partitions}.\par

Motivated by the work \citet{Gne06} in the partition context and by the product form of the EFPF of the IBP and of the Beta--Bernoulli, we consider the class of distributions for consistent exchangeable feature allocations with EFPF of the form
\begin{equation} \label{form1}
\pi _{n}(m_{1},\ldots,m_{k})=V_{n,k}\prod_{l=1}^{k}W_{m_{l}}U_{n-m_{l}},
\end{equation}
for an infinite array $V=\left( V_{n,k}:\left(n,k\right)\in\mathbb{N}\times\mathbb{N}_{0}\right) $
and two sequences $W=\left( W_{j}:j\in\mathbb{N}\right) $ and $U=\left( U_{j}:j\in\mathbb{N}_{0}\right) $
of non-negative weights, where $\mathbb{N}$ denotes the set of positive natural numbers and $\mathbb{N}_{0}=\{0\}\cup \mathbb{N}$. \par
In the feature context, we show that the IBP and the Beta--Bernoulli are the only consistent exchangeable feature allocations with form \eqref{form1}, up to randomization of their $\gamma$ and $N$ parameters respectively. Consistency and exchangeability also imply that the two sequences of weights, $W$ and $U$, must have the same form as in the IBP and in the Beta--Bernoulli, for two constants $\alpha$ and $\theta$, satisfying $\alpha \leq 1$ and $\theta > -\alpha$. In addition, $V$ must satisfy a recursion with coefficients depending on $\alpha$ and $\theta$ and the set of solutions of this recursion forms convex set. For each fixed $\alpha$ and $\theta$, we describe the extreme points of this convex set.  Their form only depends on the value of $\alpha$. For $0< \alpha <1$, the set of extreme points coincide with the family of $V$ of a 3-parameter IBP. For $\alpha =0$, this set is still continuous and coincides with the family of $V$ of the 2-parameter IBP. For $\alpha <0$ the set of extreme points is countably infinite and each extreme point corresponds to the $V$ of a Beta--Bernoulli model. To sum up, we will prove the following theorem.

\bigskip

\textsc{Theorem 1.1} 
A consistent exchangeable feature allocation has EPFP of the form \eqref{form1} iff, for some $ \alpha < 1$ and $ \theta > - \alpha $, $W_{m}=(1- \alpha)_{m-1\uparrow}$ and $U_{m}=(\theta + \alpha)_{m\uparrow}$ and the elements of $V$ satisfies the recursion 
\begin{equation} \label{form2}
V_{n,k}=\sum_{j=0}^{\infty}\binom{k+j}{j}\left(\left(\theta + \alpha \right)_{n\uparrow}\right)^{j}\left(\theta+n\right)^{k}V_{n+1,k+j}.
\end{equation}
Moreover, for fixed $\left(\alpha,\theta\right)$, the set of solutions of \eqref{form2} are 
\begin{enumerate}
\item for $0< \alpha <1$, mixtures over $\gamma$ of the $V$ of a 3-parameter IBP;
\item for $\alpha =0$, mixtures over $\gamma$ of the $V$ of a 2-parameter IBP; 
\item for $\alpha <0$, mixtures over $N$ of the $V$ of a Beta--Bernoulli model with $N$ features. 
\end{enumerate} 

\bigskip

In the next section, we prove Theorem 1.1. Specifically, in subsection 2.1, we prove the first part, characterizing $U$ and $W$ and finding the recursion for $V$. In subsection 2.2, we describe how to derive the extreme solutions of this recursion. Finally, in subsection 2.3, we study the three cases $0<\alpha<1$, $\alpha=0$, and $\alpha<0$.

\section{Proof of Theorem 1.1}
The problem is to describe all distributions for exchangeable feature allocations with EFPF \eqref{form1} subject to the consistency constraint \eqref{consistency.form}, which becomes 
\begin{equation} \label{consistency.formula}
\begin{split}
V_{n,k}  \prod_{i=1}^{k}W_{m_{l}}U_{n-m_{l}} =\sum_{j=0}^{\infty}\binom{k+j}{j} U_{n}^{j} W_{1}^{j}\sum_{\underline{z}\in\left\{ 0,1\right\} ^{k}}V_{n+1,k+j}\prod_{i=1}^{k}W_{m_{i}+z_{i}}U_{n+1-m_{i}-z_{i}},
\end{split}
\end{equation} 
for all $ n\in\mathbb{N}$, $k\in\mathbb{N}_{0}$, and for all $m_{i}\leq n$ with $i\leq k$. 
We start by noting that the representation \eqref{form1} is not unique. Specifically, we can scale
the weights in the following ways, for $\kappa>0$, and obtain the same distribution:
\begin{enumerate}
\item $\tilde{V}_{n,k}=\kappa^{-k}V_{n,k}$ and $\tilde{W}_{j}=\kappa W_{j}$;
\item $\tilde{V}_{n,k}=\kappa^{-k}V_{n,k}$ and $\tilde{U}_{j}=\kappa U_{j}$;
\item $\tilde{V}_{n,k}=\kappa^{-nk}V_{n,k}$, $\tilde{W}_{j}=\kappa^{j}W_{j}$
and $\tilde{U}_{j}=\kappa^{j}U_{j}$;
\item $\tilde{V}_{n,k}=\kappa^{-k\left(n-1\right)}V_{n,k}$, $\tilde{W}_{j}=\kappa^{j-1}W_{j}$
and $\tilde{U}_{j}=\kappa^{j}U_{j}$.
\end{enumerate}
By imposing $W_{1}=1$, we avoid the first ambiguity and with $U_{0}=1$
we fix the second one. These conditions also exclude the third ambiguity,
but do not exclude the last one, which will be fixed following Proposition 2.1.

\smallskip

\subsection{Characterization of $W$ and $U$}
The following Proposition shows that $W$ and $U$ must have a form akin to that of the IBP and $V$ is constrained to satisfy a particular recursion. In the statement of Proposition 2.1, $(x)_{n \uparrow \tau}$ denotes the generalized rising factorial.  \par

\bigskip

\textsc{Proposition 2.1.}
The weights $V$, $W$, and $U$, with the normalizations $W_{1}=U_{0}=1$,
define a consistent random feature allocation of form \eqref{form1} iff for some $a,b>0$ and $\tau\geq0$\\
(i) $W_{m}=\left(a\right)_{m-1\uparrow\tau}$, for all $m \in \mathbb{N}$;\\
(ii) $U_{m}=\left(b\right)_{m\uparrow\tau}$, for all $m \in \mathbb{N}_{0}$;\\
(iii) For all $(n,k)\in\mathbb{N}\times\mathbb{N}_{0}$, $V$ satisfies 
\begin{align} \label{form3.1}
\begin{split}
 & \sum_{j \geq 0}V_{1,j}=1 ,\\
 V_{n,k}=\sum_{j=0}^{\infty}\binom{k+j}{j} & \left(\left(b\right)_{n\uparrow\tau}\right)^{j}\left(a+b+\tau\left(n-1\right)\right)^{k}V_{n+1,k+j}.
\end{split}
\end{align}

\textsc{proof.}
The consistent exchangeable feature allocation with no features with probability one can be represented as in \eqref{form1}, with $V_{n,0}=1$ for all $n \in \mathbb{N}$ and $V_{n,k}=0$ for $k\geq 1$. 
The consistency condition \eqref{consistency.formula} for $k=1$ gives 
\begin{equation*}
V_{n,1}W_{m_{1}}U_{n-m_{1}}=\sum_{j=0}^{\infty}\left(j+1\right)V_{n+1,j+1}U_{n}^{j}\left(W_{m_{1}+1}U_{n-m_{1}}+W_{m_{1}}U_{n+1-m_{1}}\right).
\end{equation*}
This condition implies that, for all $n \in \mathbb{N}$ and for all $m_{1}\leq n$,
\begin{equation} \label{form.Prop3.1}
\frac{W_{m_{1}+1}}{W_{m_{1}}}+\frac{U_{n+1-m_{1}}}{U_{n-m_{1}}}=\frac{V_{n,1}}{\sum_{j=0}^{\infty}\left(j+1\right)V_{n+1,j+1}U_{n}^{j}}.
\end{equation}
Since the right hand side of \eqref{form.Prop3.1} does not
depend on $m_{1}$, it follows that, for all $n$ and for all $i,j\leq n$,
\begin{equation*}
\frac{W_{i+1}}{W_{i}}-\frac{W_{j+1}}{W_{j}}=\frac{U_{n+1-j}}{U_{n-j}}-\frac{U_{n+1-i}}{U_{n-i}}.
\end{equation*}
In particular, considering $n=2$, $i=2$, and $j=1$, we find
\begin{equation*}
\frac{W_{3}}{W_{2}}-W_{2}=\frac{U_{2}}{U_{1}}-U_{1}=:\tau.
\end{equation*}
For $n>1$, $i=n$, and $j=n-1$, we also obtain
\begin{equation*}
\frac{W_{n+1}}{W_{n}}-\frac{W_{n}}{W_{n-1}}=\frac{U_{2}}{U_{1}}-U_{1}=\tau,
\end{equation*}
which implies, for all $n>1$,
\begin{equation*}
\frac{W_{n+1}}{W_{n}}=\tau\left(n-1\right)+W_{2},
\end{equation*}
hence $W_{n}=\left(W_{2}\right)_{n-1\uparrow\tau}$.
In a similar manner, we consider $n>1$, $i=2$, and $j=1$ and obtain
\begin{equation*}
\frac{U_{n}}{U_{n-1}}-\frac{U_{n-1}}{U_{n-2}}=\frac{W_{3}}{W_{2}}-W_{2}=\tau.
\end{equation*}
As before, this formula implies $U_{n}=\left(U_{1}\right)_{n\uparrow\tau}$. The recursion \eqref{form3.1} follows, by rewriting \eqref{consistency.formula} as
\begin{equation*}
V_{n,k}=\sum_{j=0}^{\infty}\binom{k+j}{j}U_{n}^{j}V_{n+1,k+j}\sum_{\underline{z}\in\left\{ 0,1\right\} ^{k}}\prod_{i=1}^{k}\frac{W_{m_{i}+z_{i}}}{W_{m_{l}}}\frac{U_{n+1-m_{i}-z_{i}}}{U_{n-m_{l}}},
\end{equation*}
and by noticing that
\begin{equation*}
\sum_{\underline{z}\in\left\{ 0,1\right\} ^{k}}\prod_{i=1}^{k}\frac{W_{m_{i}+z_{i}}}{W_{m_{l}}}\frac{U_{n+1-m_{i}-z_{i}}}{U_{n-m_{l}}}=\left(U_{1}+W_{2}+\tau\left(n-1\right)\right)^{k}.
\end{equation*}
Also, $\sum_{j=0}^{\infty}V_{1,j}=1$ comes from $\sum_{j=0}^{\infty}V_{1,j}W_{1}^{j}=1$
and $W_{1}=1$. Finally, the reverse implication easily follows by checking that the probability distribution with form  \eqref{form1} and $V$, $W$ and $U$ as in the statement of the proposition satisfies the consistency condition \eqref{consistency.formula}.
\begin{flushright}
\qed
\end{flushright}

\smallskip

The last possible rescaling can now be fixed by imposing $\tau=1$. Indeed, let $W$ a sequence of weights parametrizing a feature allocation of form \eqref{form1}. From Proposition 2.1, $W_{j}=\left(a\right)_{j-1\uparrow\tau}$ for some $a>0$ and some $\tau\geq 0$. If we consider the rescaling 
$\tilde{W}_{j}=\kappa^{j-1}W_{j}$, from
$\tilde{W}_{j}=\kappa^{j-1}\left(a\right)_{j-1\uparrow\tau}=\left(\kappa a\right)_{j-1\uparrow\kappa\tau}$ and Proposition 2.1, we obtain $\tilde{W}_{j}=\left(\tilde{a}\right)_{j-1\uparrow\tilde{\tau}}$, where $\tilde{a}=\kappa a$ and $\tilde{\tau}=\kappa\tau$. Therefore, by
imposing $\tau=1$, we avoid the last ambiguity on the rescaling since $\kappa$ must be equal to $1$ and $\tilde{W}=W$. \par
We now introduce the parametrization, $\alpha:=1-a$ and $\theta:=a+b$, for $\alpha<1$ and $\theta>-\alpha$. Then, $W_{m_{l}}=\left(1-\alpha\right)_{m_{l}-1\uparrow}$ and $U_{n-m_{l}}=\left(\theta+\alpha\right)_{n-m_{l}\uparrow}$, which matches the form of the $W$ and $U$ for the IBP.

\smallskip

\subsection{General tools to derive the extreme $V$}
Let $\mathcal{V}_{\alpha,\theta}$ be the set of those elements $V \in \mathbb{R}_{+}^{\mathbb{N}\times \mathbb{N}_{0}}$ satisfying \eqref{form3.1}. Endow this set 
with the smallest $\sigma$-algebra $\mathcal{B}_{\mathcal{V}}$ that makes the maps $V \mapsto V_{n,k}$ measurable and define the barycenter $V^{\mu}$ of each measure $\mu$ on $\mathcal{B}_{\mathcal{V}}$ as the pointwise average,
\begin{equation}
V_{n,k}^{\mu}=\int_{\mathcal{V}_{\alpha,\theta}} \!\! V_{n,k}\,\mu(\mathrm{d}V).
\end{equation}
It is easy to check that $\mathcal{V}_{\alpha,\theta}$ is a convex set, i.e., for all probability measures $\mu$ on $\mathcal{B}_{\mathcal{V}}$, $V^{\mu}\in \mathcal{V}_{\alpha,\theta}$ (see Appendix \ref{App.Fact1}). The goal of this section is to check that this set is also a simplex and to describe its extreme elements.  \par

Given a measurable space of functions with the convex structure just defined, Dynkin \citet{Dyn78} describes a general theory which can be applied to show that this set is also a simplex and to determine its extreme points. Similar results have been studied or rediscovered by many works, see references in \citet{Gne06}.  To apply the results of \citet{Dyn78} to our problem, we will follow the same strategy used by \citet{Gne06}. Rather than studying $\mathcal{V}_{\alpha,\theta}$ directly, we consider another space, isomorphic to $\mathcal{V}_{\alpha,\theta}$ and easier to study, and we find its extreme points applying the results by \citet{Dyn78}. \par

Let $(\mathbb{N}_{0}^{\infty},\mathcal{C}(\mathbb{N}_{0}^{\infty}))$ be the infinite product space of $\mathbb{N}_{0}$, endowed with its cylinder $\sigma$-algebra. To each $V \in \mathcal{V}_{\alpha,\theta}$ we associate a Markov law, $P_{V}$, on this space. Specifically, writing $K_{n}:\mathbb{N}_{0}^{\infty} \rightarrow \mathbb{N}_{0}$ for the $n$-th coordinate projection on the product space, the Markov law associated to $V$ has the initial distribution given by
\begin{equation} \label{form_init.dist}
P_{V}\left(K_{1}=j\right)=V_{1,j},
\end{equation}
and transition probabilities
\begin{equation} \label{form_trans.prob}
P_{V}\left(K_{n+1}=j+k|K_{n}=k\right)=\binom{k+j}{j}\left(\left(\alpha + \theta \right)_{n\uparrow}\right)^{j}\left(\theta + n\right)^{k}\frac{V_{n+1,k+j}}{V_{n,k}},
\end{equation}
if $j\geq0$ and $0$ otherwise. Let $\mathcal{P}_{\mathcal{V}_{\alpha,\theta}}=\{P_{V}:V\in \mathcal{V}_{\alpha,\theta}\}$ be the set of Markov laws. The map $T:\mathcal{V}_{\alpha,\theta} \rightarrow \mathcal{P}_{\mathcal{V}_{\alpha,\theta}}$, defined by $T(V)=P_{V}$ is a convex isomorphism (see Appendix \ref{App.Fact1} for a proof). Hence, if $P_{V}$ is extreme in $\mathcal{P}_{\mathcal{V}_{\alpha,\theta}}$, so is $V$ in $\mathcal{V}_{\alpha,\theta}$. We now describe how to find the extreme points of $\mathcal{P}_{\mathcal{V}_{\alpha,\theta}}$. Before that, we remark that, given an EFPF with form \eqref{form1} parametrized by $V$, it is straightforward to show that $K_{n}$ corresponds to the number of features in the corresponding random feature allocation of $n$ individuals, i.e., $K_{n}$ is the cardinality of $\bigcup_{1 \leq i \leq n} X_{i}$. \par

As we will see from Proposition 2.2, 
for every $n\in \mathbb{N}$, $\mathcal{F}_{n}=\sigma (K_{n},K_{n+1},\ldots)$ is a sufficient $\sigma$-algebra for  $\mathcal{P}_{\mathcal{V}_{\alpha,\theta}}$. Hence, for each $n\in \mathbb{N}$, there exists a common $\mathcal{F}_{n}$-measurable regular conditional probability $Q_{n}:\mathbb{N}_{0}^{\infty} \times \mathcal{C}(\mathbb{N}_{0}^{\infty})\rightarrow [0,1]$ for $\mathcal{P}_{\mathcal{V}_{\alpha,\theta}}$ given $\mathcal{F}_{n}$, such that, for all $P_{V} \in \mathcal{P}_{\mathcal{V}_{\alpha,\theta}}$  
and $A \in \mathcal{C}(\mathbb{N}_{0}^{\infty})$,
\begin{equation} \label{Q_n}
Q_{n}(\omega,A)=P_{V}((K_{m})_{m \in \mathbb{N}} \in A | \mathcal{F}_{n})(\omega),
\end{equation}
for $P_{V}$-almost all $\omega \in \mathbb{N}_{0}^{\infty}$. 
In order to avoid having to repeat uninteresting measure-theoretic details, 
when $A' \in \sigma(K_1,\dotsc,K_n)$, we will take advantage of the Markov property 
of $(K_m)_{m \in \mathbb{N}}$
to assume that 
\begin{equation}
Q_n(\omega,A') 
= P((K_{m})_{m \in \mathbb{N}} \in A' | \mathcal{F}_{n})(\omega)
= P((K_{m})_{m \in \mathbb{N}} \in A' | K_n)(\omega)
\end{equation}
for all $\omega \in \mathbb{N}_{0}^{\infty}$,
where we have dropped the $V$ from the notation $P_V$ 
in order to highlight the independence of the cotransition probabilities under $P_V$ from $V$ itself.
This is justified because the equality holds for all $P \in \mathcal{P}_{\mathcal{V}_{\alpha,\theta}}$.

Associated to each Markov kernel $Q_n$, there is a Markov operator $\Pi_n$ given by
\begin{equation} \label{Markov.operators}
\Pi _{n} f(\omega)=\int f(\omega') Q_{n}(\omega,\mathrm{d}\omega'),
\end{equation}
for all  $f$ bounded $\mathcal{F}_{n}$-measurable real functions. Henceforth, for every $\sigma$-algebra $\mathcal{F}$, we will simply write  $f\in \mathcal{F}$ to denote that $f$ is bounded and $\mathcal{F}$-measurable.  The sequence $(\mathcal{F}_{n},\Pi _{n})_{n \in \mathbb{N}}$ forms a \emph{specification} in $(\mathbb{N}_{0}^{\infty},\mathcal{C}(\mathbb{N}_{0}^{\infty}))$ (see Appendix \ref{App.Fact1} for a proof). We can apply Theorem 5.1 of \citet{Dyn78}, which states that  $(\Pi _{n})_{n \in \mathbb{N}}$ is an asymptotically H-sufficient statistic, which in turn means (see also Section 4.4 of \citet{Dyn78}) that, for all $P_{V}$ that are extreme, 
\begin{equation} \label{form.Dynkyn}
P_{V}(\{\omega \in \mathbb{N}_{0}^{\infty}: \forall f \in \mathcal{C}(\mathbb{N}_{0}^{\infty}), \; \underset{n\rightarrow \infty}{\lim} \Pi _{n} f(\omega)=\textstyle\int \!f \, \mathrm{d}P_{V}\}) = 1.
\end{equation}
A path $\omega \in \mathbb{N}_{0}^{\infty}$ induces a Markov law $P_{V}\in\mathcal{P}_{\mathcal{V}_{\alpha,\theta}}$ and is said to be \emph{regular} iff 
for all $f \in \mathcal{C}(\mathbb{N}_{0}^{\infty})$, $\underset{n\rightarrow \infty}{\lim} \Pi _{n} f(\omega)=\int \! f \, \mathrm{d}P_{V}$. 
The set of point in $\mathcal{P}_{\mathcal{V}_{\alpha,\theta}}$ that are induced by regular paths is called the \emph{maximal boundary}. The set of extreme points, also called the \emph{minimal boundary}, is the subset of the maximal boundary, corresponding to those points $P_{V}$ that also satisfy \eqref{form.Dynkyn}, i.e., they assign probability 1 to the set of regular paths inducing them. \par

In our context, to identify the maximal boundary, it is enough to check \eqref{form.Dynkyn} for all functions $f\in\mathcal{C}(\mathbb{N}_{0}^{\infty})$ 
that are indicators of cylinder sets of the form $K_{n}^{-1}\{k\}$ for $n\in\mathbb{N}$ and $k\in\mathbb{N}_{0}$. That is, the elements belonging to the maximal boundary are those $P_{\bar{V}} \in \mathcal{P}_{\mathcal{V}_{\alpha,\theta}}$ such that, for some $\omega \in \mathbb{N}_{0}^{\infty}$,
\begin{equation*}
\underset{m\rightarrow \infty}{\lim} P(K_{n}=k | \mathcal{F}_{m})(\omega)=\underset{m\rightarrow \infty}{\lim} P(K_{n}=k | K_{m})(\omega)=P_{\bar{V}}(K_{n}=k),
\end{equation*} %
for all $(n,k) \in \mathbb{N}\times\mathbb{N}_{0}$. 
To find the extremes measures of $\mathcal{P}_{\mathcal{V}_{\alpha,\theta}}$, we compute the cotransition (backwards) probabilities of $(K_{n})_{n \in \mathbb{N}}$. 

\bigskip

\textsc{Proposition 2.2.} 
The cotransition probabilities are 
\begin{equation}
P(K_{n}=k|K_{m}=l)=\frac{d_{n,k}^{m,l}}{d^{m,l}} d^{n,k},
\end{equation}
for $n<m$ and $k\leq l$, while the distribution of $K_{n}$ under $P_{V} \in \mathcal{P}_{\mathcal{V}_{\alpha,\theta}}$ is
\begin{equation}
P_{V}(K_{n}=k)=V_{n,k} d^{n,k},
\end{equation}
where 
\begin{equation}
d^{m,l}=( (\theta + 1)_{m-1\uparrow} + \sum_{j=1}^{m-1}(\theta + \alpha)_{m-j\uparrow}(\theta+1+m-j)_{j-1\uparrow} )^{l},
\end{equation}
and
\begin{equation}
d_{n,k}^{m,l}= \binom{l}{k} ((\theta+n)_{m-n\uparrow})^{k}(\sum_{j=1}^{m-n}(\alpha + \theta)_{m-j\uparrow}(\theta+m-1)_{j-1\downarrow})^{l-k}.
\end{equation}

\medskip
\textsc{proof.}
See Appendix \ref{App.prop3.2}.
\begin{flushright}
\qed
\end{flushright} \par
\smallskip
\noindent Note that the cotransition probabilities are independent of $V$.

\smallskip

\subsection{Characterization of $V_{n,k}$}
In this section, we study the three cases $0<\alpha<1$, $\alpha=0$, and $\alpha<0$ separately.
Recall that a path $\omega=(\omega_{1},\omega_{2},\ldots) \in \mathbb{N}_{0}^{\infty}$ is regular and induces $\bar{V}\in \mathcal{V}_{\alpha,\theta}$ if the limit
\begin{equation}\label{thelimit} 
\underset{m\rightarrow \infty}{\lim} P(K_{n}=k|K_{m}=\omega_{m})=\underset{m\rightarrow \infty}{\lim}  \frac{d_{n,k}^{m,\omega_{m}}}{d^{m,\omega_{m}}}  d^{n,k} = \bar{V}_{n,k} d^{n,k}
\end{equation}
exists for all $(n,k)$. In this case, $P_{\bar{V}}$ belongs to the maximal boundary of $\mathcal{P}_{\mathcal{V}_{\alpha,\theta}}$. If $P_{\bar{V}}$ also assigns probability one to the set of regular paths inducing it, then $P_{\bar{V}}$ is extreme. 

\subsubsection{Case $0<\alpha<1$}
For $(\alpha,\theta)$ fixed, s.t. $0<\alpha<1$ and $\theta>-\alpha$, let $V^{3IBP,\alpha,\theta}(\gamma)$ be the $V$ of the 3-parameter IBP, defined as
\begin{align*}
V^{3IBP,\alpha,\theta}_{n,k}(\gamma)& = \frac{1}{k!}\Big(\frac{\gamma}{\left(\theta+1\right)_{n-1\uparrow}}\Big)^{k}\exp\Big(-\sum_{i=1}^{n}\gamma\frac{\left(\alpha+\theta\right)_{i-1\uparrow}}{\left(1+\theta\right)_{i-1\uparrow}}\Big) \\
& = \frac{1}{k!} \Big(\frac{\gamma}{(\theta + 1)_{n-1\uparrow}}\Big)^{k} \exp \Big(-\gamma \Big( \frac{\Gamma (\theta +1) \Gamma( \alpha+\theta+n)}{\alpha  \Gamma( \alpha+\theta) \Gamma( \theta +n)} - \frac{\theta}{\alpha} \Big) \Big),
\end{align*}
for all $\gamma\geq 0$. 
Define  $\mathcal{P}_{V^{3IBP,\alpha,\theta}}=\{P_{V^{3IBP,\alpha,\theta}(\gamma)} \in \mathcal{P}_{\mathcal{V}_{\alpha,\theta}}:\gamma \geq 0\}$. \par

\bigskip

\textsc{Proposition 2.3.}
Let $0<\alpha<1$ and $\theta>-\alpha$.
\begin{enumerate}
\item[a)] The elements of the set $\mathcal{P}_{V^{3IBP,\alpha,\theta}}$ belong to the maximal boundary of $\mathcal{P}_{\mathcal{V}_{\alpha,\theta}}$ and they are induced by those paths $ w \in \mathbb{N}^{ \infty}_{0} $ s.t. $ \frac{w_{m}}{m^{\alpha}} \rightarrow c$, where $c=\frac{\gamma \Gamma(\theta +1)}{\alpha \Gamma(\alpha + \theta)}$; 
\item[b)] The elements of $\mathcal{P}_{V^{3IBP,\alpha,\theta}}$ also belong to the minimal boundary of $\mathcal{P}_{\mathcal{V}_{\alpha,\theta}}$, i.e., they are extreme points of $\mathcal{P}_{\mathcal{V}_{\alpha,\theta}}$;
\item[c)] The elements of $\mathcal{P}_{V^{3IBP,\alpha,\theta}}$ are the only extreme points, i.e., $\mathcal{P}_{V^{3IBP,\alpha,\theta}}$ coincides with the maximal and the minimal boundary.
\end{enumerate} \par
\smallskip
\textsc{proof.}
\begin{enumerate}
\item[a)]  In Appendix \ref{App.prop3.3.1}, we check that 
\begin{equation} \label{limit0<alpha<1}
\lim_{\substack{m\to\infty \\ \frac{\omega_{m}}{m^{\alpha}} \to c}} \frac{d_{n,k}^{m,\omega_{m}}}{d^{m,\omega_{m}}} = V^{3IBP,\alpha,\theta}_{n,k} \Big( \frac{c\alpha\Gamma(\alpha+\theta)}{\Gamma(\theta +1)} \Big).
\end{equation}
\item[b)] From Theorem 4 of \citet{Ber15}, it follows that
\begin{equation*}
P_{V^{3IBP,\alpha,\theta} ( \frac{c\alpha\Gamma(\alpha+\theta)}{\Gamma(\theta +1)} )} \Big( \frac{K_{n}}{n^{\alpha}}\rightarrow c \Big)=1.
\end{equation*}
\item[c)] In Appendix \ref{App.prop3.3.1}, we show that the elements of  $\mathcal{P}_{V^{3IBP,\alpha,\theta}}$ are the only ones belonging to the maximal boundary, i.e., there are not any other regular paths except those of part (a).
\end{enumerate} 
\begin{flushright}
\qed
\end{flushright}

In Proposition 2.3, the case $\gamma=0$ corresponds to the degenerate feature allocation with no features with probability one, corresponding to $V_{n,0}=1$ and $V_{n,k}=0$ for all $n\in \mathbb{N}$ and $k\geq 1$. This solution is induced by the path $\omega_{m}=0$ for all $m \in \mathbb{N}$, which has probability one under this degenerate law. 

\subsubsection{Case $\alpha=0$}
For $\theta$ fixed and positive, the $V$ of the 2-parameter IBP are of the form
\begin{align*}
V^{2IBP,\theta}_{n,k}(\gamma) & = \frac{1}{k!}\bigg(\frac{\gamma}{\left(\theta+1\right)_{n-1\uparrow}}\bigg)^{k}\exp\bigg(-\sum_{i=1}^{n}\gamma\frac{\left(\theta\right)_{i-1\uparrow}}{\left(1+\theta\right)_{i-1\uparrow}}\bigg) \\
& = \frac{1}{k!}\bigg(\frac{\gamma}{(\theta+1)_{n-1\uparrow}}\bigg)^{k} \exp\bigg(-\gamma\sum_{i=1}^{n}\frac{\theta}{\theta+i-1}\bigg),
\end{align*}
with the convention that, when $\gamma=0$, we recover the degenerate feature allocation with no features. 
Define $\mathcal{P}_{V^{2IBP,\theta}}=\{P_{V^{2IBP,\theta}(\gamma)} \in \mathcal{P}_{\mathcal{V}_{0,\theta}}:\gamma \geq 0\}$. \par

\bigskip

\textsc{Proposition 2.4.}
Let $\alpha =0$ and $\theta>0$. 
\begin{enumerate}
\item[a)] The elements of the set $\mathcal{P}_{V^{2IBP,\theta}}$ belong the maximal boundary of $\mathcal{P}_{\mathcal{V}_{0,\theta}}$ and they are induced by paths $ w \in \mathbb{N}^{ \infty}_{0} $ s.t. $\frac{w_{m}}{\log(m)} \rightarrow \gamma $; 
\item[b)] The elements of $\mathcal{P}_{V^{2IBP,\theta}}$ also belong to the minimal boundary of $\mathcal{P}_{\mathcal{V}_{0,\theta}}$, i.e., they are extreme points of $\mathcal{P}_{\mathcal{V}_{0,\theta}}$;
\item[c)] The elements of $\mathcal{P}_{V^{2IBP,\theta}}$ are the only extreme points, i.e., $\mathcal{P}_{V^{3IBP,\theta}}$ coincides with the maximal and the minimal boundary.
\end{enumerate} \par
\smallskip
\textsc{proof.}
\begin{enumerate}
\item[a)] In Appendix \ref{App.prop3.3.2}, we check that 
\begin{equation} \label{limit.alpha=0}
\lim_{\substack{m\to\infty \\ \frac{\omega_{m}}{\log(m)} \to \gamma}} \frac{d_{n,k}^{m,\omega_{m}}}{d^{m,\omega_{m}}} = V^{2IBP,\theta}_{n,k} ( \gamma ).
\end{equation}
\item[b)] This also follows from Theorem 4 of \citet{Ber15}, which establish that
\begin{equation*}
P_{V^{2IBP,\theta} ( \gamma )} \Big( \frac{K_{n}}{\log(n)}\rightarrow \gamma \Big)=1.
\end{equation*}
\item[c)] In Appendix \ref{App.prop3.3.2}, we check that there are not any other regular paths but those of part (a).
\end{enumerate}
\begin{flushright}
\qed
\end{flushright}

\subsubsection{Case $\alpha <0$}
From formula \eqref{Beta--Bernoulli}, we see that the Beta--Bernoulli is of form \eqref{form1}, with $V$ of the form
\begin{equation}
V_{n,k}^{BB,\alpha,\theta}(N)=\frac{\binom{N}{k} \big(\frac{-\alpha \Gamma(\theta+\alpha)}{\Gamma(\theta+\alpha+n)}\big)^k}{\big(\frac{\Gamma(\theta+\alpha)\Gamma(\theta+n)}{\Gamma(\theta+\alpha+n)\Gamma(\theta)}\big)^{N}}.
\end{equation}
for all $N\in \mathbb{N}$. As before, when $N=0$, we consider the feature allocation with no feature with probability one. 
Define $\mathcal{P}_{V^{BB,\alpha,\theta}}=\{ P_{V^{BB,\alpha,\theta}(N)} \in \mathcal{P}_{\mathcal{V}_{\alpha,\theta}}: N \in \mathbb{N}_{0} \}$.

\bigskip

\textsc{Proposition 2.5.}
Let $\alpha <0$ and $\theta > -\alpha$. 
\begin{enumerate}
\item[a)] The elements of the set $\mathcal{P}_{V^{BB,\alpha,\theta}}$ belong the maximal boundary of $\mathcal{P}_{\mathcal{V}_{\alpha,\theta}}$  and they are induced by paths $ w \in \mathbb{N}^{ \infty}_{0}$ s.t. $w_{m} \rightarrow N$; 
\item[b)] The elements of $\mathcal{P}_{V^{BB,\alpha,\theta}}$ also belong to the minimal boundary of $\mathcal{P}_{\mathcal{V}_{\alpha,\theta}}$, i.e., they are extreme points of $\mathcal{P}_{\mathcal{V}_{\alpha,\theta}}$;
\item[c)] The elements of $\mathcal{P}_{V^{BB,\alpha,\theta}}$ are the only extreme points, i.e., $\mathcal{P}_{V^{BB,\alpha,\theta}}$ coincides with the maximal and the minimal boundary.
\end{enumerate} \par
\smallskip
\textsc{proof.}
\begin{enumerate}
\item[a)] In Appendix \ref{App.prop3.3.3}, we check that 
\begin{equation} \label{limit.alpha<0}
\lim_{\substack{m\to\infty \\ \omega_m \to N}} \frac{d_{n,k}^{m,\omega_{m}}}{d^{m,\omega_{m}}} = V^{BB,\alpha,\theta}_{n,k} (N).
\end{equation}
\item[b)] This follows since under a Beta--Bernoulli model with $N$ features, $K_{m}\rightarrow N$ a.s. Indeed, the probability of each feature, $q_{j}$, is a.s. strictly positive, being Beta distributed. The probability of this feature having all zeros in a $m$-individuals allocation is $(1-q_{k})^{m}$, which tends to zero as $m\rightarrow \infty$. 
\item[c)] In Appendix \ref{App.prop3.3.3}, we check that there are not any other regular paths but those of part (a).
\end{enumerate} 
\begin{flushright}
\qed
\end{flushright}

\section{Discussion}
In this work, we have considered the class of consistent exchangeable feature allocations with EPFP of the form \eqref{form1}. While this is a tractable family, the only elements of this class are mixtures over $\gamma$ of the 2 and 3-parameter IBP or mixtures over $N$ of the Beta--Bernoulli model. From both an applied and theoretical perspective, it would be of interest to have larger but still tractable classes of exchangeable feature allocations. Finding new tractable priors for feature models is still an active area of research. A possible direction of research would be to study a more general class than \eqref{form1}, with form $V_{n,k}\prod_{l=1}^{k}W_{n,m_{l}}$, for a triangular array $W=\left( W_{n,k}:n\in\mathbb{N},0\leq k \leq n\right)$. However, a characterization of $W$ in this case would seem to be much more complicated than Proposition 2.1.

\appendix

\section{Appendix: Some Facts and Proofs}

\subsection{EFPF of the Beta--Bernoulli model} \label{App.BetaBernoulli}
The \emph{Beta--Bernoulli} is described by considering a finite space of features, which can numbered using the integers in $[N]$, where $N$ is the cardinality of the feature set. To each feature we associate a random number $q_{j}$ with distribution $\text{Beta}(\eta_{1},\eta_{2})$. Each individual $X_{i}$ possesses feature $j$ with probability $q_{j}$. Let $Z_{i,j}$ be a binary random variable denoting the presence or absence of feature $j$ in individual $i$. Then, the Beta--Bernoulli model can be written as
\begin{align*}
Z_{i,j}|q_{j} & \sim \text{Bernoulli}(p_{j}) \;  & i=1, & \ldots,n;  \; j=1,\dots,N; \\
q_{j}|\eta_{1},\eta_{2} & \sim \text{Beta}(\eta_{1},\eta_{2})   & & j=1,\dots,N.
\end{align*}
The conditional probability that $Z=(Z_{i,j})_{i\leq n,j\leq N}$ is equal to $z=(z_{i,j})_{i\leq n,j\leq N}$ given $q=\left(q_{1},\ldots,q_{N}\right)$ is
\begin{equation*}
p(Z=z|q)=\prod_{j=1}^{N} \prod_{i=1}^{n} \text{Bernoulli}(z_{i,j}|q_{j}).
\end{equation*}
Integrating $q$ out, we obtain the probability mass function of $Z$
\begin{equation*}
p(Z=z)=\left( \frac{\Gamma(\eta_{1}+\eta_{2})}{\Gamma(\eta_{1})\Gamma(\eta_{2})} \right) ^{N} \prod_{i=1}^{N}\frac{\Gamma(m_{i}+ \eta_{1})\Gamma(n-m_{i}+\eta_{2})}{\Gamma(n+\eta_{1} + \eta_{2})},
\end{equation*} 
where $m_{i}=\sum_{j=1}^{n}z_{i,j}$. If $(m_{1},\ldots,m_{N})$ has $k$ non-zero entries, this probability becomes
\begin{equation*}
\Big( \frac{\Gamma(\eta_{1}+\eta_{2})}{\Gamma(\eta_{1})\Gamma(\eta_{2})\Gamma(n+\eta_{1} + \eta_{2}) } \Big) ^{N} \left(\Gamma(\eta_{1})\Gamma(n+\eta_{2})\right)^{N-k} \prod_{i=1}^{k}\Gamma(m_{i}+ \eta_{1})\Gamma(n-m_{i}+\eta_{2}).
\end{equation*}
Finally, taking into account all $\binom{N}{k}$ possible uniform orderings of the $N-k$ features not possessed by any individual, which give rise to the same uniformly ordered feature allocation, we obtain
\begin{equation*}
\binom{N}{k} \left( \frac{\Gamma(\eta_{1}+\eta_{2})}{\Gamma(\eta_{1})\Gamma(\eta_{2})\Gamma(n+\eta_{1} + \eta_{2}) } \right) ^{N} \left(\Gamma(\eta_{1})\Gamma(n+\eta_{2})\right)^{N-k} 
\cdot \prod_{i=1}^{k}\Gamma(m_{i}+ \eta_{1})\Gamma(n-m_{i}+\eta_{2}),
\end{equation*}
which can be be rewritten as in formula \eqref{Beta--Bernoulli}, by using rising factorials and by changing the parametrization to $\alpha=-\eta_{1}$ and $\theta=\eta_{2}+\eta_{1}$, with $\alpha<0$ and $\theta>-\alpha$. \par

\medskip

\subsection{Some Facts} \label{App.Fact1} 

\textsc{Proposition A.1.}
 $\mathcal{V}_{\alpha,\theta}$ is a convex set. \par
\medskip
\textsc{proof.}
We want to show that $\mathcal{V}_{\alpha,\theta}$ is a convex set, i.e., for all probability measures $\mu$ on $\mathcal{B}_{\mathcal{V}}$, $V^{\mu}\in \mathcal{V}_{\alpha,\theta}$. We have
\begin{align*}
V_{n,k}^{\mu} & = \int_{\mathcal{V}_{\alpha,\theta}}V_{n,k}\mu(\mathrm{d}V) \\
 & =\int_{\mathcal{V}_{\alpha,\theta}}\sum_{j=0}^{ \infty} \binom{k+j}{j} ((\alpha + \theta)_{n \uparrow})^{j}  (\theta+n)^{k}  V_{n+1,k+j}\mu(\mathrm{d}V) \\
 & =\sum_{j=0}^{ \infty} \binom{k+j}{j} ((\alpha + \theta)_{n \uparrow})^{j}  (\theta+n)^{k}  \int_{\mathcal{V}_{\alpha,\theta}} V_{n+1,k+j}\mu(\mathrm{d}V) \\
 & =\sum_{j=0}^{ \infty} \binom{k+j}{j} ((\alpha + \theta)_{n \uparrow})^{j} (\theta+n)^{k}  V_{n+1,k+j}^{\mu},
\end{align*} 
for all $(n,k)$, where the first and last equality follow from the definition of barycenter, and the second from the monotone convergence theorem. In a similar manner,
\begin{align*}
\sum_{j=0}^{\infty}V_{1,j}^{\mu} & =  \sum_{j=0}^{\infty} \int_{\mathcal{V}_{\alpha,\theta}}V_{1,j}\mu(\mathrm{d}V) \\
 &  = \int_{\mathcal{V}_{\alpha,\theta}} \sum_{j=0}^{\infty} V_{1,j}\mu(\mathrm{d}V) \\
 &  = \int_{\mathcal{V}_{\alpha,\theta}} 1 \mu(\mathrm{d}V) = 1.
\end{align*} 
\begin{flushright}
\qed
\end{flushright}

\bigskip

\indent \textsc{Proposition A.2.}
$T(V)=P_{V}$ is an isomorphism between convex sets. \par
\smallskip
 \textsc{proof.}
According to \cite[][pg~706]{Dyn78}, the map $T(V)=P_{V}$ is a convex isomorphism if $T$ is invertible and $T$ and $T^{-1}$ are measurable and preserve the convex structure. $T$ is 1-1 from Proposition 2.2 and it is onto by construction. We prove $T$ is measurable and preserves the convex structure, proving that the same is true for $T^{-1}$ can be done in similar way. \par
$\mathcal{P}_{\mathcal{V}_{\alpha,\theta}}$ is endowed with the smallest $\sigma$-algebra $\mathcal{B}_{\mathcal{P}}$ that makes the maps $P_{V}\mapsto P_{V}(A)$ measurable for all $A \in \mathcal{C}(\mathbb{N}_{0}^{\infty})$. A generator of this $\sigma$-algebra is composed by sets $\{ P_{V} \in \mathcal{P}_{\mathcal{V}_{\alpha,\theta}} : P_{V} (K_{n} =k) \leq x \}$ for $(n,k)\in \mathbb{N}\times\mathbb{N}_{0}$ and $x\in[0,1]$. The inverse image under $T$ of this set is $ \{ V \in \mathcal{V}_{\alpha,\theta}:V_{n,k}d_{n,k}\leq x \}$ (see Proposition 2.2 for the definition of $d_{n,k}$), which lies in $\mathcal{B}_{\mathcal{V}}$. Hence, $T$ is measurable. \par
$T$ preserves the convex structure if, for every measure $\mu$ on $\mathcal{B}_{\mathcal{V}}$, we have $T(V^{\mu})=P^{\mu'}$, where $\mu '$ the push-forward measure of  $\mu$  on $\mathcal{B}_{\mathcal{P}}$ (i.e., $\mu '=\mu \circ T^{-1}$), and $P^{\mu'}$ is the barycenter of $\mu'$, defined as
\begin{equation}
P^{\mu'}(A)=\int_{\mathcal{P}_{\mathcal{V}_{\alpha,\theta}}} P(A)\mu'(\mathrm{d}P),
\end{equation}
for all $ A \in \mathcal{C}(\mathbb{N}_{0}^{\infty})$. Using the change of variable formula, it is easy to check that $T$ preserves the convex structure. Indeed, considering cylinder sets of the form $K_{n}^{-1}\{k\}$ for $n\in\mathbb{N}$ and $k\in\mathbb{N}_{0}$, we have 
\begin{align*}
P^{\mu '}(K_{n}=k) & = \int_{\mathcal{P}_{\mathcal{V}_{\alpha,\theta}}} P(K_{n}=k)\mu'(\mathrm{d}P)= \int_{\mathcal{P}_{\mathcal{V}_{\alpha,\theta}}} P(K_{n}=k)\mu \circ T^{-1}(\mathrm{d}P) \\
   & =  \int_{\mathcal{V}_{\alpha,\theta}} d_{n,k}V_{n,k} \mu(\mathrm{d}V)=d_{n,k}V^{\mu}_{n,k}.
\end{align*}
Hence, $T(V^{\mu})=P^{\mu'}$.
\begin{flushright}
\qed
\end{flushright}

\bigskip

\textsc{Proposition A.3.}  
$(\mathcal{F}_{n},\Pi _{n})_{n \in \mathbb{N}}$  forms a specification in $(\mathbb{N}_{0}^{\infty},\mathcal{C}(\mathbb{N}_{0}^{\infty}))$. \par
\smallskip
\textsc{proof.}
According to \citet{Dyn78}, section 5.1, given a directed set $L$ and a measurable space $(\Lambda,\mathcal{F})$, a \emph{specification} on this is space $(\mathcal{F}_{\Lambda},\Pi_{\Lambda})_{\Lambda \in L}$ is a family of sub-$\sigma$-algebras and Markov operators satisfying
\begin{enumerate}
\item[(i)] $\mathcal{F}_{\Lambda'}\subseteq \mathcal{F}_{\Lambda}$, if $\Lambda' \succeq \Lambda$;
\item[(ii)] $\Pi_{\Lambda'}\Pi_{\Lambda} = \Pi_{\Lambda'}$, if $\Lambda'\succeq \Lambda$;
\item[(iii)] $\Pi_{\Lambda}f \in \mathcal{F}_{\Lambda}$, for all $f \in \mathcal{F}$;
\item[(iv)] $\Pi_{\Lambda}f = f$, for all $f \in \mathcal{F}_{\Lambda}$.
\end{enumerate}
In our context, with $L=\mathbb{N}$ and the sub-$\sigma$-algebras and Markov operators defined in section 3.1, formula \eqref{Markov.operators}, (i), (ii), and (iv) follow immediately. To check (ii), it is enough to check for indicators of measurable sets. In particular, for $f=\mathbbm{1}_{A}$, with $A \in \mathcal{C}(\mathbb{N}_{0}^{\infty})$, we must check 
\begin{equation*}
\int_{\mathbb{N}_{0}^{\infty}}Q_{n}(\omega ',A)Q_{n+1}(\omega,\mathrm{d}\omega ')=Q_{n+1}(\omega,A).
\end{equation*}
Indeed, it is enough to check this condition for a thin cylinder $A$ of the form $K_{1}^{-1}\{k_{1}\}\cap K_{2}^{-1}\{k_{2}\}\cap \ldots K_{m}^{-1}\{k_{m}\}$ for $m>n+1$ and $k_{i} \in \mathbb{N}_{0}$ for all $i\leq m$, \par
\small
\begin{align*}
 \int_{\mathbb{N}_{0}^{\infty}} & Q_{n}(\omega ',A)Q_{n+1}(\omega,\mathrm{d}\omega ')  \\
& = \int_{ \mathbb{N}_{0}^{\infty}} P(K_{1}=k_{1},\ldots,K_{m}=k_{m}|\mathcal{F}_{n})(\omega ') P((K_{l})_{l \in \mathbb{N}} \in \mathrm{d}\omega '|\mathcal{F}_{n+1})(\omega) \\
&  = \int _{ \mathbb{N}_{0}^{\infty}} P(K_{1}=k_{1},\ldots,K_{n-1}=k_{n-1}|K_{n}=k_{n}) \\
& \;\;\;\;\;\;\;\;\;\;\;\;\;\;\;\;\;\;\;\;\;\;  \cdot\mathbbm{1}(\omega_{n}'=k_{n})\ldots\mathbbm{1}(\omega_{m}'=k_{m}) P((K_{l})_{l \in \mathbb{N}} \in \mathrm{d}\omega '|\mathcal{F}_{n+1})(\omega) \\
& = P(K_{1}=k_{1},\ldots,K_{n-1}=k_{n-1}|K_{n}=k_{n})  \\
& \;\;\;\;\;\;\;\;\;\;\;\;\;\;\;\;\;\;\;\;\;\;  \cdot\int _{ \mathbb{N}_{0}^{\infty}}\mathbbm{1}(\omega_{n}'=k_{n})\ldots\mathbbm{1}(\omega_{m}'=k_{m}) P((K_{l})_{l \in \mathbb{N}} \in \mathrm{d}\omega '|\mathcal{F}_{n+1})(\omega) \\
& =  P(K_{1}=k_{1},\ldots,K_{n-1}=k_{n-1}|K_{n}=k_{n}) P(K_{n}=k_{n}|\mathcal{F}_{n+1})(\omega) \\
& \;\;\;\;\;\;\;\;\;\;\;\;\;\;\;\;\;\;\;\;\;\;  \cdot\mathbbm{1}(\omega_{n+1}=k_{n+1})\ldots\mathbbm{1}(\omega_{m}=k_{m})\\
& =  P(K_{1}=k_{1},\dots,K_{m}=k_{m}| \mathcal{F}_{n+1})(\omega)=Q_{n+1}(\omega,A).
\end{align*} \par
\normalsize
\begin{flushright}
\qed
\end{flushright}

\medskip

\subsection{Proof of Proposition 2.2} \label{App.prop3.2}
First, note that for $m>n$ and $l\geq k$  $P_{V}(K_{m}=l|K_{n}=k)=V_{m,l} d_{n,k}^{m,l}$, for a function $d_{n,k}^{m,l}$ independent of $V$. Indeed, from \eqref{form_trans.prob}, the probability of a path $(k_{n+1},k_{n+2},\ldots,k_{m-1},l)$ depends only on the last $V_{m,l}$. Summing over all possible paths from $K_{n}=k$ to $K_{m}=l$, we see that $P_{V}(K_{m}=l|K_{n}=k)$ must be of the form $V_{m,l} d_{n,k}^{m,l}$. In addition, by considering $P_{V}(K_{m}=l)=\sum_{i=0}^{l} P_{V}(K_{m}=l|K_{1}=i)\cdot P_{V}(K_{1}=i)$, $P_{V}(K_{m}=l)$ must be of the form $V_{m,l} d^{m,l}$. Also, from
\begin{equation*}
P_{V}\left(K_{m}=l\right)=\sum_{j=0}^{l}P_{V}\left(K_{m}=l|K_{m-1}=j\right)\cdot P_{V}\left(K_{m-1}=j\right),
\end{equation*} 
it follows that, for $l>2$, the function $d^{m,l}$ must satisfy
\begin{equation*}
V_{m,l}d^{m,l}=\sum_{j=0}^{l}\frac{V_{m,l}}{V_{m-1,j}}\tbinom{l}{l-j}(\left(\alpha + \theta \right)_{m-1\uparrow})^{l-j}\left(\theta+m-1\right)^{j}V_{m-1,j}d^{m-1,j},
\end{equation*}
which gives the following the recursion
\begin{equation*} 
d^{m,l}=\sum_{j=0}^{l}\tbinom{l}{l-j}(\left( \alpha + \theta \right)_{m-1\uparrow})^{l-j}\left( \theta +m-1\right)^{j}d^{m-1,j}.
\end{equation*}
Substituting $d^{m-1,j}$, we find \par
\small
\begin{equation*}
d^{m,l}=\sum_{j=0}^{l}\tbinom{l}{l-j}(\left( \alpha + \theta \right)_{m-1\uparrow})^{l-j}\left( \theta + m-1\right)^{j}\sum_{i=0}^{j}\binom{j}{j-i}(\left( \alpha + \theta \right)_{m-2\uparrow})^{j-i}\left(\theta+m-2\right)^{i}d^{m-2,i}.
\end{equation*} \par
\normalsize
\noindent Grouping together all coefficient multiplying $d^{m-2,k}$ on the right hand side $\left(0\leq k\leq l\right)$, we find %
\small
\begin{align*}
d_{n-2,k}^{m,l}& =  \tbinom{l}{l-k}\left(\left(\alpha + \theta \right)_{m-1\uparrow}\right)^{l-k}\left(\theta +m-1\right)^{k}\left( \theta+m-2\right)^{k} \\
 &  +\tbinom{l}{l-k-1}\left(\left(\alpha + \theta \right)_{m-1\uparrow}\right)^{l-k-1}\left(\theta +m-1\right)^{k+1}\left(\left( \alpha + \theta \right)_{m-2\uparrow}\right)\left(\theta +m-2\right)^{k} \\
 &  +\tbinom{l}{l-k-2}\left(\left(\alpha + \theta \right)_{m-1\uparrow}\right)^{l-k-2}\left(\theta +m-1\right)^{k+2}\binom{k+2}{2}\left(\left(\alpha + \theta \right)_{m-2\uparrow}\right)^{2}\left(\theta +m-2\right)^{k} \\
 &  \qquad \vdots \\
 &  +\left(\left(\alpha+\theta \right)_{m-1\uparrow}\right)\left(\theta+m-1\right)^{l-1}\binom{l-1}{l-1-k}\left(\left( \alpha + \theta \right)_{m-2\uparrow}\right)^{l-1-k}\left(\theta +m-2\right)^{k} \\
 &   +\left(\theta+m-1\right)^{l}\binom{l}{l-k}\left(\left(\alpha+\theta \right)_{m-2\uparrow}\right)^{l-k}\left(\theta+m-2\right)^{k} \\
 & = \left(\left(\theta+m-1\right)\left(\theta+m-2\right)\right)^{k}\left(\left(\alpha +\theta \right)_{m-1\uparrow}+\left(\theta +m-1\right)\left(\alpha + \theta \right)_{m-2\uparrow}\right)^{l-k}\binom{l}{l-k}.
\end{align*} \par
\normalsize
\noindent So, the recursion for $d^{m,l}$ becomes %
\small
\begin{equation*}
d^{m,l}=\sum_{j=0}^{l}\binom{l}{l-j}\left(\left(\theta+m-1\right)\left(\theta+m-2\right)\right)^{j}\left(\left(\alpha + \theta \right)_{m-1\uparrow}+\left(\theta+m-1\right)\left(\alpha +\theta \right)_{m-2\uparrow}\right)^{l-j} d^{m-2,j}.
\end{equation*} \par
\normalsize
\noindent In the same manner, we find %
\small
\begin{align*}
d_{m-3,k}^{m,l}= & \left(\left(\theta+m-1\right)\left(\theta+m-2\right)\left(\theta+m-3\right)\right)^{k}\binom{l}{l-k} \\
& \cdot\left(\left(\alpha + \theta \right)_{m-1\uparrow}+\left(\theta+m-1\right)\left(\alpha+\theta\right)_{m-2\uparrow}+\left(\theta+m-1\right)\left(\theta+m-2\right)\left(\alpha+\theta\right)_{m-3\uparrow}\right)^{l-k}.
\end{align*} \par
\normalsize
\noindent Finally, we obtain,
\begin{align*}
d_{n,k}^{m,l} & =\binom{l}{l-k}((\theta+m-1)_{m-n\downarrow})^{k}(\sum_{j=1}^{m-n}(\alpha + \theta )_{m-j\uparrow}(\theta+m-1)_{j-1\downarrow})^{l-k} \\
&=\binom{l}{k}((\theta+n)_{m-n\uparrow})^{k}(\sum_{j=1}^{m-n}(\alpha + \theta )_{m-j\uparrow}(\theta+m-1)_{j-1\downarrow})^{l-k}.
\end{align*}
In addition,
\begin{align*}
d^{m,l}= & \sum_{i=0}^{l} d_{1,i}^{m,l} = \sum_{i=0}^{l}\binom{l}{i}((\theta+1)_{m-1\uparrow})^{i}(\sum_{j=1}^{m-1}(\alpha + \theta )_{m-j\uparrow}(\theta+m-1)_{j-1\downarrow})^{l-i} \\
 & =( (\theta +1)_{m-1\uparrow} + \sum_{j=1}^{m-1}(\alpha + \theta)_{m-j\uparrow}(\theta+1+m-j)_{j-1\uparrow} )^{l}.
\end{align*}
\begin{flushright}
\qed
\end{flushright}
\medskip

\subsection{Proof of Proposition 2.3} \label{App.prop3.3.1}

We begin with some technical results about asymptotic equivalence of functions:
Write $f \approx g$ to denote that $f(m) / g(m) \rightarrow 1$ as $m \to \infty$,
and note that $f_i \approx g_i$ implies $f_1+f_2 \approx f_1 + g_2 \approx g_1 + g_2$ and $f_1 f_2 \approx f_1 g_2 \approx g_1 g_2$.
In general, it does not hold that $f_1^{f_2} \approx g_1^{g_2}$.
The following two results characterize special cases: \par
\medskip

\textsc{Lemma A.4.}\label{expeqv}
Let $g(m) \to \infty$ as $m \to \infty$, let $f \approx g$ and $h/g \to c$ for some constant $c \ge 0$.
Then, for every $p,q \in \mathbb{R}$, we have
\begin{equation}
\left ( \frac {f(m)-p}{f(m)-q} \right )^{h(m)}
\approx
\left ( \frac {g(m)-p}{g(m)-q} \right )^{c\, g(m)}
\to e^{c\,(q-p)}
\end{equation}
\smallskip

\textsc{proof}
We prove only the first equivalence, because the limiting exponential form is well known.
Taking logarithms, we have
\begin{align}
&\log \left [ \left ( \frac {f(m)-p}{f(m)-q} \right )^{h(m)}
\left ( \frac {g(m)-q} {g(m)-p} \right )^{c\, g(m)}
\right ] \label{lograt}
\\&=
c\, g(m)  \, \log \frac {(f(m)-p)\,(g(m)-q)}{(f(m)-q)\,(g(m)-p)} 
\\&\qquad + (h(m) - c\, g(m)) \, \log \frac {f(m)-p}{f(m)-q}.  
\end{align}
The arguments to the logarithms can we written as
\begin{align}
\frac {(f(m)-p)\,(g(m)-q)}{(f(m)-q)\,(g(m)-p)} 
&= 1 + \frac {  (f(m)- g(m))\,(p-q)} {(f(m)-q)(g(m)-p)} \label{eqa}
\end{align}
and
\begin{align}
\frac {f(m)-p}{f(m)-q} 
= 1 + \frac {q-p}{f(m)-q}. \label{eqb}
\end{align}
Using the fact that
$z(m) \to 0$ implies $\log (1+z(m)) \approx z(m)$,
and that both terms \eqref{eqa} and \eqref{eqb} 
converge to one, 
it follows that
\begin{align}
\eqref{lograt} \approx
 c\, g(m) \frac {  (f(m)- g(m))\,(p-q)} {(f(m)-q)(g(m)-p)} 
 + (h(m)-c\, g(m)) \, \frac {q-p}{f(m)-q}.
\end{align}
It is straightforward to show that the r.h.s.\ converges to 0.
\begin{flushright}
\qed
\end{flushright} \par
\smallskip

\textsc{lemma A.5.}\label{lemdom}
Let $f(m) \to \infty$ as $m \to \infty$, let $g(m)/f(m) \to \infty$ as $m \to \infty$, and let $h \approx f$.
For every $p > q$, 
\begin{equation*}
\biggl ( \frac {g(m)}{f(m)} \biggr )^k \biggl ( \frac {h(m)-p}{h(m) - q} \biggr)^{g(m)} \to 0  \qquad \text{as $m \to \infty$.}
\end{equation*} \par
\smallskip
\textsc{proof.}
Taking logarithms
\begin{align*}
&k \log \frac {g(m)}{f(m)} + g(m) \log \Bigl \{ 1 + \frac {q-p}{h(m)-q} \Bigr \}
\\&\approx k \log \frac {g(m)}{f(m)} + (q-p) \frac {g(m)}{h(m)-q}
\\&\approx k \log \frac {g(m)}{f(m)} + (q-p) \frac {g(m)}{f(m)} \to -\infty
\end{align*}
as $m \to \infty$, completing the proof.
\begin{flushright}
\qed
\end{flushright} \par

\bigskip

We now proceed to prove each part of Proposition 2.3: \par
\smallskip
a)  We must check the limit \eqref{limit0<alpha<1}. From Proposition 2.2, 
\small
\begin{align}
\frac{d_{n,k}^{m,\omega_{m}}}{d^{m,\omega_{m}}} 
&= \binom{\omega_{m}}{k} \frac{\left[\left(\theta+m-1\right)_{m-n\downarrow}\right]^{k}\left[\sum_{i=1}^{m-n}\left(\alpha+\theta\right)_{m-i\uparrow}\left(\theta+m-1\right)_{i-1\downarrow}\right]^{\omega_{m}-k}}{\left[\left(\theta+m-1\right)_{m-1\downarrow}+\sum_{i=1}^{m-1}\left(\alpha+\theta\right)_{m-i\uparrow}\left(\theta+m-1\right)_{i-1\downarrow}\right]^{\omega_{m}}} \nonumber\\
& =  \binom{\omega_{m}}{k}\left[\frac{\left(\theta+m-1\right)_{m-n\downarrow}}{\sum_{i=1}^{m-n}\left(\alpha+\theta\right)_{m-i\uparrow}\left(\theta+m-1\right)_{i-1\downarrow}}\right]^{k}  \nonumber\\
& \qquad \cdot \left[\frac{\sum_{i=1}^{m-n}\left(\alpha+\theta\right)_{m-i\uparrow}\left(\theta+m-1\right)_{i-1\downarrow}}{\left(\theta+m-1\right)_{m-1\downarrow}+\sum_{i=1}^{m-1}\left(\alpha+\theta\right)_{m-i\uparrow}\left(\theta+m-1\right)_{i-1\downarrow}}\right]^{\omega_{m}} \nonumber\\
\begin{split}
& =  
\binom{\omega_{m}}{k}\left[\frac{\frac{\Gamma\left(\theta+m\right)}{\Gamma\left(\theta+n\right)}}{\frac{\Gamma\left(\alpha+\theta+m\right)}{\alpha\cdot\Gamma\left(\alpha+\theta\right)}-\frac{\Gamma\left(\theta+m\right)\cdot\Gamma\left(\alpha+\theta+n\right)}{\alpha\cdot\Gamma\left(\alpha+\theta\right)\cdot\Gamma\left(\theta+n\right)}}\right]^{k} \\
& \qquad  \cdot \left[\frac{\frac{\Gamma\left(\alpha+\theta+m\right)}{\alpha\cdot\Gamma\left(\alpha+\theta\right)}-\frac{\Gamma\left(\theta+m\right)\cdot\Gamma\left(\alpha+\theta+n\right)}{\alpha\cdot\Gamma\left(\alpha+\theta\right)\cdot\Gamma\left(\theta+n\right)}   }{\frac{\Gamma\left(\theta+m\right)}{\Gamma\left(\theta+1\right)}+\frac{\Gamma\left(\alpha+\theta+m\right)}{\alpha\cdot\Gamma\left(\alpha+\theta\right)}-\frac{\Gamma\left(\theta+m\right)\cdot\Gamma\left(\alpha+\theta+1\right)}{\alpha\cdot\Gamma\left(\alpha+\theta\right)\cdot\Gamma\left(\theta+1\right)}}\right]^{\omega_{m}}, \label{aoeu1}
\end{split}
\end{align} 
\normalsize 
where the third equality follows from the identity
\small
\begin{equation*}
\sum_{i=1}^{m-n}\left(\alpha+\theta\right)_{m-i\uparrow}\left(\theta+m-1\right)_{i-1\downarrow}=\frac{\Gamma\left(\alpha+\theta+m\right)}{\alpha\Gamma\left(\alpha+\theta\right)}-\frac{\Gamma\left(\theta+m\right)\Gamma\left(\alpha+\theta+n\right)}{\alpha\Gamma\left(\alpha+\theta\right)\Gamma\left(\theta+n\right)},
\end{equation*} 
\normalsize
which itself arises from rewriting the sum as a difference of two infinite hypergeometric series evaluated at 1 
and applying the Gauss theorem for hypergeometric series.

\indent
Using the asymptotic equivalence 
$\Gamma(m+\delta) \approx \Gamma(m)m^{\delta}$
and limit $(m+\theta)^\alpha - m^\alpha \to 0$,
the Stirling formula $\binom{\omega_{m}}{k} \approx  \frac{1}{k!}\omega_{m}^{k}$ for the binomial coefficient,
and then the limit $\frac{\omega_{m}}{m^{\alpha}}\to c$,
the first line of \eqref{aoeu1} can be simplified to yield a limiting form:
\begin{align}
&\binom{\omega_{m}}{k}
   \Biggl[\frac{\frac{\Gamma\left(\theta+m\right)}{\Gamma\left(\theta+n\right)}}
                     {\frac{\Gamma\left(\alpha+\theta+m\right)}{\alpha\cdot\Gamma\left(\alpha+\theta\right)}
                          -\frac{\Gamma\left(\theta+m\right)\cdot\Gamma\left(\alpha+\theta+n\right)}
                                  {\alpha\cdot\Gamma\left(\alpha+\theta\right)\cdot\Gamma\left(\theta+n\right)}}
    \Biggr]^{k} 
\nonumber
\\&\approx
\binom{\omega_{m}}{k} \Biggl [ \frac{\frac{1}{\Gamma(\theta+n)}}{\frac{m^{\alpha}}{\alpha\Gamma(\alpha+\theta)}-\frac{\Gamma(\alpha+\theta+n)}{\alpha\Gamma(\alpha+\theta)\Gamma(\theta+n)}} \Biggr ] ^{k} 
\nonumber
\\&\approx
 \frac{1}{k!} \Biggr[ \omega_{m}m^{-\alpha} \cdot \frac{\frac{1}{\Gamma(\theta+n)}}{\frac{1}{\alpha\Gamma(\alpha+\theta)}-\frac{\Gamma(\alpha+\theta+n)}{m^{\alpha}\alpha\Gamma(\alpha+\theta)\Gamma(\theta+n)}}  \Biggr] ^{k} 
\label{eqv3}
\\&\approx
 \frac{1}{k!}  \Biggr [ c \cdot \frac{\frac{1}{\Gamma(\theta+n)}}{\frac{1}{\alpha\Gamma(\alpha+\theta)}-\frac{\Gamma(\alpha+\theta+n)}{m^{\alpha}\alpha\Gamma(\alpha+\theta)\Gamma(\theta+n)}}  \Biggr ] ^{k} 
\rightarrow \frac{1}{k!} \left( c \frac{\alpha\Gamma(\alpha+\theta)}{\Gamma(\theta+n)}\right)^{k}.
\nonumber
\end{align}
Similarly, the second line of \eqref{aoeu1}
can be simplified
by Lemma~\ref{expeqv}, 
using the asymptotic equivalence
$\Gamma(m+\delta) \approx \Gamma(m)m^{\delta}$
and the limits $(m+\theta)^\alpha - m^\alpha \to 0$ %
and $\frac{\omega_{m}}{m^{\alpha}}\to c$, 
to yield
\begin{align}
&\Biggl[\frac{\frac{\Gamma\left(\alpha+\theta+m\right)}{\alpha\cdot\Gamma\left(\alpha+\theta\right)}-\frac{\Gamma\left(\theta+m\right)\cdot\Gamma\left(\alpha+\theta+n\right)}{\alpha\cdot\Gamma\left(\alpha+\theta\right)\cdot\Gamma\left(\theta+n\right)}   }{\frac{\Gamma\left(\theta+m\right)}{\Gamma\left(\theta+1\right)}+\frac{\Gamma\left(\alpha+\theta+m\right)}{\alpha\cdot\Gamma\left(\alpha+\theta\right)}-\frac{\Gamma\left(\theta+m\right)\cdot\Gamma\left(\alpha+\theta+1\right)}{\alpha\cdot\Gamma\left(\alpha+\theta\right)\cdot\Gamma\left(\theta+1\right)}}\Biggr]^{\omega_{m}}
\label{secondeqv}
\\&\approx
  \Biggl[ \frac{m^{\alpha}-\frac{\Gamma(\alpha+\theta+n)}{\Gamma(\theta+n)}}{m^{\alpha} - \frac{\theta\Gamma(\alpha+\theta)}{\Gamma(\theta+1)} }  \Biggr]^{cm^{\alpha}}
 \to 
\exp \left \{ c \left( \frac{\theta\Gamma(\alpha+\theta)}{\Gamma(\theta+1)}-\frac{\Gamma(\alpha+\theta+n)}{\Gamma(\theta+n)} \right) \right \}.
\nonumber
\end{align}
Substituting back into \eqref{aoeu1}, we obtain the $V$ of the 3-parameter IBP.
\smallskip

c) To check that the only regular paths are those paths $\omega \in \mathbb{N}_{0}^{ \mathbb{N}}$ such that $\frac{w_{m}}{m^{\alpha}} \rightarrow c$ for some $c\geq 0$, 
suppose otherwise; i.e., let $\omega \in \mathbb{N}_0^{ \mathbb{N}}$ be a regular path, but assume $\frac{w_{m}}{m^{\alpha}}$ does not converge to some finite $c\geq 0$. 
If $(\frac{w_{m}}{m^{\alpha}})_{m\in \mathbb{N}}$ has at least two distinct subsequential limits, then, from the proof of part (a), 
we see that $d_{n,k}^{m,\omega_{m}}/d^{m,\omega_{m}}$ has at least two distinct subsequential limits,
a contradiction, and so $\frac{w_{m}}{m^{\alpha}} \rightarrow \infty$. 
But then it follows from equations \eqref{aoeu1}, \eqref{eqv3}, and \eqref{secondeqv};
the asymptotic equivalence
$\Gamma(m+\delta) \approx \Gamma(m)m^{\delta}$
and limit $(m+\theta)^\alpha - m^\alpha \to 0$;
and finally an application of Lemma~\ref{lemdom} that 
$\frac{d_{n,k}^{m,\omega_{m}}}{d^{m,\omega_{m}}} \to 0$ as $m \to \infty$ for every $k \in \mathbb{N}_{0}$.  As these limits must define a probability distribution, this is a contradiction, completing the proof.
\qed

\subsection{Proof of Proposition 2.4} \label{App.prop3.3.2}
a) We must check the limit \eqref{limit.alpha=0}. From Proposition 2.2,
\small
\begin{align*}
\frac{d_{n,k}^{m,\omega_{m}}}{d^{m,\omega_{m}}}&=\frac{\binom{\omega_{m}}{k}\left[\left(\theta+m-1\right)_{m-n\downarrow}\right]^{k}\left[\sum_{i=1}^{m-n}\left(\theta\right)_{m-i\uparrow}\left(\theta+m-1\right)_{i-1\downarrow}\right]^{l-k}}{\left[\left(\theta+m-1\right)_{m-1\downarrow}+\sum_{i=1}^{m-1}\left(\theta\right)_{m-i\uparrow}\left(\theta+m-1\right)_{i-1\downarrow}\right]^{\omega_{m}}} \\
& =\binom{\omega_{m}}{k}\Bigg[\frac{\frac{\Gamma\left(\theta+m\right)}{\Gamma\left(\theta+n\right)}}{\frac{\Gamma\left(\theta+m\right)}{\Gamma\left(\theta\right)}\sum_{i=1}^{m-n}\frac{1}{\theta+m-i}}\Bigg]^{k}\Bigg[\frac{\frac{\Gamma\left(\theta+m\right)}{\Gamma\left(\theta\right)}\sum_{i=1}^{m-n}\frac{1}{\theta+m-i}}{\frac{\Gamma\left(\theta+m\right)}{\Gamma\left(\theta\right)}\sum_{i=1}^{m-1}\frac{1}{\theta+m-i}+\frac{\Gamma\left(\theta+m\right)}{\Gamma\left(\theta+1\right)}}\Bigg]^{\omega_{m}},
\end{align*}
\normalsize 
where the second equality follows from the identity
\begin{equation*}
\sum_{i=1}^{m-n}\left(\theta\right)_{m-i\uparrow}\left(\theta+m-1\right)_{i-1\downarrow}=\left(\theta\right)_{m\uparrow}\sum_{i=1}^{m-n}\frac{1}{\theta+m-i}.
\end{equation*}
Using the Stirling formula for the binomial coefficient, $\binom{\omega_{m}}{k} \approx  \frac{1}{k!}\omega_{m}^{k}$, and the identity
$\sum_{i=1}^{m}\frac{1}{\theta+m-i} = \sum_{i=1}^{m-n}\frac{1}{\theta+m-i} + \sum_{i=1}^{n}\frac{1}{\theta+n-i}$,
we have
\begin{align*}
\frac{d_{n,k}^{m,\omega_{m}}}{d^{m,\omega_{m}}}
&\approx\frac{1}{k!}\Bigg[\omega_{m}\frac{\frac{\Gamma\left(\theta\right)}{\Gamma\left(\theta+n\right)}}{\sum_{i=1}^{m-n}\frac{1}{\theta+m-i}}\Bigg]^{k}\Bigg[\frac{\sum_{i=1}^{m-n}\frac{1}{\theta+m-i}}{\sum_{i=1}^{m}\frac{1}{\theta+m-i}}\Bigg]^{\omega_{m}}
\\&\approx\frac{1}{k!}\Bigg[\frac{\omega_{m}}{\sum_{i=1}^{m-n}\frac{1}{\theta+m-i}} \frac{\Gamma\left(\theta\right)}{\Gamma\left(\theta+n\right)} \Bigg]^{k}\Bigg[\frac{\sum_{i=1}^{m-n}\frac{1}{\theta+m-i}}{\sum_{i=1}^{m-n}\frac{1}{\theta+m-i} +  \sum_{i=1}^{n}\frac{1}{\theta+n-i}}\Bigg]^{\omega_{m}}.
\end{align*}
Therefore, by Lemma~\ref{expeqv} and the fact that 
$\log(m) \approx \sum_{i=1}^{m-n}\frac{1}{\theta+m-i}$
and $\omega_{m}/ \log(m) \to \gamma $, 
we have
\begin{align*}
\frac{d_{n,k}^{m,\omega_{m}}}{d^{m,\omega_{m}}}
&\approx\frac{1}{k!}\left[\gamma  
           \frac{\Gamma\left(\theta\right)}{\Gamma\left(\theta+n\right)} \right]^{k}
             \left[\frac {\log(m)}{\log(m) + \sum_{i=1}^{n}\frac{1}{\theta+n-i} } \right]^{\gamma \log(m)}
\\&\to
\frac{1}{k!}\left[\gamma \frac{\Gamma\left(\theta\right)}{\Gamma\left(\theta+n\right)} \right]^{k} 
           \exp \Big( - \gamma\sum_{i=1}^{n}\frac{1}{\theta+n-i} \Big),
\end{align*}
as $m\to\infty$, recovering the $V$ of the 2-parameter IBP.\par
\smallskip
c) To check that the only regular paths are those paths $\omega \in  \mathbb{N}^ \mathbb{N}_{0}$ such that $\frac{\omega_{m}}{\log(m)} \rightarrow \gamma$ for some $\gamma\geq 0$, we can repeat the same argument as in the proof of Proposition 2.3, part (c). First, we note that if $(\frac{w_{m}}{\log(m)})_{m\in\mathbb{N}}$ has at least two distinct subsequential limits, then $\omega$ cannot be regular, because the proof of point (a) shows that there will be two distinct induced laws, a contradiction. If $\frac{\omega_{m}}{\log(m)}\rightarrow \infty$ as $m\rightarrow\infty$, then it follows again from Lemma~\ref{lemdom} that
$\frac{d_{n,k}^{m,\omega_{m}}}{d^{m,\omega_{m}}} \to 0$ as $m \to \infty$ for all $k\in \mathbb{N}_{0}$, a contradiction.

\begin{flushright}
\qed
\end{flushright}

\subsection{Proof of Proposition 2.5} \label{App.prop3.3.3}

a) We must check the limit \eqref{limit.alpha<0}. Starting with Proposition 2.2 and following similar steps as for the case $0<\alpha<1$, we obtain the approximation
\begin{equation} \label{proof.alpha<0}
\frac{d_{n,k}^{m,\omega_{m}}}{d^{m,\omega_{m}}}  \approx  \binom{\omega_{m}}{k} \Bigg( \frac{\frac{1}{\Gamma(\theta+n)}}{\frac{m^{\alpha}}{\alpha\Gamma(\alpha+\theta)}-\frac{\Gamma(\alpha+\theta+n)}{\alpha\Gamma(\alpha+\theta)\Gamma(\theta+n)}} \Bigg) ^{k}  \Bigg( \frac{m^{\alpha}-\frac{\Gamma(\alpha+\theta+n)}{\Gamma(\theta+n)}}{m^{\alpha} - \frac{\theta\Gamma(\alpha+\theta)}{\Gamma(\theta+1)} }  \Bigg)^{\omega_{m}},
\end{equation}
assuming $\omega_{m} \rightarrow N$ and $\alpha<0$.  Taking the limit as $m \rightarrow \infty$, we obtain
\begin{equation*} 
\frac{d_{n,k}^{m,\omega_{m}}}{d^{m,\omega_{m}}} \rightarrow\binom{N}{k}\left[\frac{-\alpha\Gamma\left(\alpha+\theta\right)}{\Gamma\left(\alpha+\theta+n\right)}\right]^{k} 
 \left[\frac{\Gamma\left(\alpha+\theta+n\right)\Gamma\left(\theta\right)}{\Gamma\left(\alpha+\theta\right)\Gamma\left(\theta+n\right)}\right]^{N}.
\end{equation*}\par

\smallskip
c) By the a.s.\ monotonicity of regular paths, as $m \to \infty$, the number of features $\omega_m$ either diverges or converges to a finite (integer) limit. The divergent paths cannot be regular for $\alpha<0$, because for these paths, \eqref{proof.alpha<0} diverges as $m\rightarrow \infty$. Hence, the only regular paths are those of part (a).  

\begin{flushright}
\qed
\end{flushright}

\section*{Acknowledgements}

S. Favaro is supported by the European Research Council through StG N-BNP 306406. D.M. Roy is funded by an NSERC Discovery Grant, Connaught New Researcher Award, AFOSR grant FA9550-15-1-0074, and Newton Alumni grant. Y.W. Teh research leading to these results has received funding from the European Research Council under the European Union's Seventh Framework Programme (FP7/2007-2013) ERC grant agreement no. 617071.

\end{document}